\documentclass[reqno,12pt]{amsart} 
\usepackage{amsmath,hyperref}
\usepackage{tikz}
\usepackage{graphicx} 
\newcommand{\Cross}{\mathbin{\tikz [x=1.8ex,y=1.8ex,line width=.2ex]
\draw (0,0) -- (1,1) (0,1) -- (1,0);}}

\newtheorem{theorem}{Theorem}
\newtheorem{proposition}{Proposition}

\theoremstyle{remark}

\newtheorem{example}{Example}

\numberwithin{equation}{section}

\newcommand{\la}{\lambda}
\newcommand{\ga}{\gamma}
\newcommand{\C}{\mathbb C}
\newcommand{\N}{\mathbb N}

\newcommand{\abs}[1]{\lvert#1\rvert}
\newcommand\cpgn[3]{{C}_{#1,#2}(#3)}
\newcommand\tcgn[2]{{T}_{#1}(#2)}
\newcommand\cps[2]{\mathcal{C}_{#1,#2}}
\newcommand\tcs[1]{\mathcal{T}_{#1}}
\newcommand{\beqs}{\begin{equation*}}
\newcommand{\eeqs}{\end{equation*}}
\newcommand{\beq}{\begin{equation}}
\newcommand{\eeq}{\end{equation}}
\newcommand\eqn[1]{(\ref{#1})}
\newcommand\twidit[1]{\overset {\text{\lower 3pt\hbox{$\sim$}}}#1}

\allowdisplaybreaks

\author{Frank Garvan}
\address{Department of Mathematics, University of Florida,
358 Little Hall, PO box 118105, Gainesville, FL 32611-8105, USA}
\email{fgarvan@ufl.edu}
\urladdr{http://www.qseries.org/fgarvan/}
\thanks{Partly supported by a grant from the Simon's Foundation (\#318714).}

\author{Michael J.\ Schlosser}
\address{Fakult\"at f\"ur Mathematik, Universit\"at Wien,
Oskar-Morgenstern-Platz~1, A-1090 Vienna, Austria}
\email{michael.schlosser@univie.ac.at}
\urladdr{http://www.mat.univie.ac.at/{\textasciitilde}schlosse}
\thanks{Partly supported by FWF Austrian Science Fund
grant F50-08 within the SFB
``Algorithmic and enumerative combinatorics''.}

\title[Interpretations of Ramanujan's tau function]
{Combinatorial interpretations of\\ Ramanujan's tau function}



\begin{document}

\begin{abstract}
We use a $q$-series identity by Ramanujan to give a combinatorial
interpretation of Ramanujan's tau function which involves $t$-cores
and a new class of partitions which we call $(m,k)$-capsids.
The same method can be applied in conjunction with other related identities
yielding alternative combinatorial interpretations
of the tau function. 
\end{abstract}

\maketitle

\section{Introduction}\label{secintro}

Ramanujan's~\cite{R1} tau function $\tau(n)$ is defined to be
the coefficient of $q^n$ in the $24$th power of Dedekind's eta function
$\eta(q):=q^{\frac 1{24}}(q;q)_\infty$, i.e.,
\beqs
\sum_{n\ge 1}\tau(n)q^n=q(q;q)^{24}_\infty.
\eeqs
Here we are using standard $q$-series notation (cf.\ \cite{GR})
\beqs
(a;q)_n:=\prod_{j=0}^{n-1}(1-aq^j), 
\eeqs
where $a\in\C$, $0<|q|<1$ and $n$ is an nonnegative integer
(for $n=0$ the expression $(a;q)_0$ reduces to an empty product
which is defined to be $1$) or infinity. For brevity, we also
frequently use the compact notation
\beqs
(a_1,\dots,a_m;q)_n:=(a_1;q)_n\dots(a_m;q)_n. 
\eeqs

The tau function possesses very nice arithmetic properties,
see \cite{Ra}. In particular, $\tau(n)$ is a multiplicative function,
as originally observed by Ramanujan and later proved by Mordell~\cite{M}.
Lehmer~\cite{L} conjectured that the tau function never vanishes,
which is still open. The tau function oscillates and assumes
positive and negative integer values. Dyson~\cite[Eq.(2), p.636]{D} 
revealed a nice
formula for $\tau(n)$ which is actually a special case of one of
the Macdonald identities~\cite{Ma}.
Although the latter formula does have a strong combinatorial flavour,
as it involves an alternating multiple sum of rational numbers
the formula is not easily seen to be related to a direct counting
of special combinatorial objects.

We are interested in finding combinatorial interpretations for $\tau(n)$.
We succeed by employing $q$-series identities of a specific type,
together with the use
of generating functions for $t$-cores and for a new class of
partitions which we call $(m,k)$-capsids.
The main idea is to split the generating function of $\tau(n)$
using an identity of the type \eqn{AB}
in a nonnegative component and a nonpositive component,
see \eqn{eq:tauid},
where all components admit easy combinatorial
interpretations in terms of $t$-cores and $(m,k)$-capsids.

\section{Partitions, $t$-cores, and $(m,k)$-capsids}
\label{sectptns}

Let $\N$ denote the set of positive integers.
A {\em partition} $\la$ of a positive integer $n$ (cf.\ \cite{A0}) is
a nonincreasing sequence of positive integers (called the {\em parts}),
$\la=(\la_1,\la_2,\dots,\la_l)$, such that $\la_1+\la_2+\dots+\la_l=n$.
If $\la$ is a partition of $n$, then we write $\la\vdash n$ or $|\la|=n$.
We can identitfy a partition $\la=(\la_1,\la_2,\dots,\la_l)$
with its {\em Ferrers diagram}, denoted by $[\la]$, defined to be
the set of cells $\{(i,j)\in\N^2: 1\le i\le l,\,1\le j\le \la_i\}$.
For a cell $(i,j)\in [\la]$,
the {\em hook} $H_\la(i,j)$ of $\la$ is the following subset of the
Ferrers diagram $[\la]$:
$H_\la(i,j)=\{(i,t)\in [\la]:t\ge j\}\cup
\{(s,j)\in [\la]:s>i\}$. The cardinality of the hook $H_\la(i,j)$
is the {\em hooklength} $h_\la(i,j)$.

Sometimes it is convenient to write a partition $\la$ not as a
nonincreasing sequence
of positive integers but rather as a increasing sequence of integers
with finitely many nonzero multiplicities, formally written in
the form $\la=(1^{\mu_1},2^{\mu_2},\dots)$.
Here, each positive integer $j\in\N$ has an nonnegative
multiplicity $\mu_j$, correpsonding to the number of times
the part $j$ occurs in $\la$.
If $\la\vdash n$, we must have $\sum_{j\ge 1}j\mu_j=n$.

For our combinatorial interpretations of Ramanujan's tau function
we will use two special classes of partitions: $t$-cores and
$(m,k)$-capsids\footnote{In medical microbiology, the {\em core} is the
inner structure of a viral particle (which consists of structural proteins
and genom), such as a rotavirus; it is surrounded by an (inner and an outer)
{\em capsid}, the protein shell of a virus.},
the latter are introduced here for the first time.

$t$-Cores have been introduced in \cite[Sec.~2.7]{JK}
where they have been shown to be useful for the recursive
evaluation of the irreducible charcters of the symmetric group.
A partition $\la$ is a {\em $t$-core} if, and only if, $\la$
contains no hooks whose lengths are multiples of $t$.
Let $c(t,n)$ denote the number of $t$-core partitions of $n$.
The generating function for $t$-cores was computed in \cite{GKS}:
\beqs
T_t(q)=1+\sum_{n=1}^\infty c(t,n)q^n=\frac{(q^t;q^t)_\infty^t}{(q;q)_\infty}.
\eeqs
See also \cite[Theorem 2.7.17, p.80]{JK}.
$t$-Cores have attracted broad interest and were studied from 
various points of view, in particular, with focus on some of their 
combinatorial \cite{G1} or analytic \cite{An,GO} properties. 

Let $m,k\in\N$, $m\ge 2$.
For $0< k < m$ we say that a partition $\pi$ is an {\em $(m,k)$-capsid}
if, and only if, the possible parts of $\pi$ are $m-k$ or are congruent
to $0$ or $k$ mod $m$, and satisfy the following two conditions:

\smallskip
(i) if $\mu_{m-k}=0$, i.e., $m-k$ is not a part,
then all parts are congruent to $k$ mod $m$

\smallskip
(ii) if $\mu_{m-k}>0$, then $m-k$ is the smallest part
and the largest part congruent to $0$ mod $m$ is
$\le m\cdot\mu_{m-k}$ and all parts congruent to $k$ mod $m$
(different from $m-k$, if $k=m/2$) are $> m\cdot\mu_{m-k}$.

\medskip
Let $\gamma(m,k,n)$ be the number of $(m,k)$-capsid partitions of $n$.
Then, from the combinatorial definition it is clear that we have the
following generating function
\begin{equation}\label{capsid}
{C}_{m,k} (q)= 1 + \sum_{n=1}^\infty \gamma(m,k,n) q^n
=
\sum_{n=0}^\infty \frac{q^{(m-k)n}}
{(q^m;q^m)_n (q^{mn+k};q^m)_\infty}.
\end{equation}

\begin{example}\label{ex:51}
Let $m=5$ and $k=1$. 
There are seven $(5,1)$-capsid partitions of $16$:
\beqs
(1^{16}), \quad (4^{4}), \quad (1^{10},6),\quad (1^{4},6^{2}),\quad 
(4,6^2),\quad
(1^5,11),\quad (16),
\eeqs
so that $a(5,1,16)=7$.                                                         
\end{example}

It turns out the generating function ${C}_{m,k}(q)$ is also an 
infinite product. 
By \cite[Eq.~(4.2)]{Be-Ga08b} 
we have
\beq\label{Cgenf}
\frac{ (at;q)_\infty }{ (a;q)_\infty (t;q)_\infty }
=
\sum_{n=0}^\infty \frac{t^n}{(q)_n (aq^n;q)_\infty},
\eeq
which is a generalization of \cite[Eq.~(5.7)]{Ga88b}. This identity
also follows immediately from the $q$-binomial
theorem~\cite[Eq.~(II.3)]{GR}. Now the
$(q,a,t)\mapsto (q^m,q^k,q^{m-k})$ case of
\eqn{Cgenf}, together with \eqn{capsid},
immediately establishes the following result:
\begin{proposition}\label{prop:caps}
The generating function for $(m,k)$-capsids is
\beqs
{C}_{m,k}(q)=\frac{ (q^m;q^m)_\infty }
{ (q^k;q^m)_\infty (q^{m-k};q^m)_\infty }.
\eeqs
\end{proposition}
Since ${C}_{m,k}(q)={C}_{m,m-k}(q)$, by inspection,
we immediately deduce that the numbers of $(m,k)$- and
$(m,m-k)$-capsid partitions of $n$ are equinumerous:
\beq\label{symcaps}
\ga(m,k,n)=\ga(m,m-k,n).
\eeq
This easy consequence of Proposition~\ref{prop:caps}
appears to be not at all combinatorially obvious.
In Section \ref{sectsyms} we give a bijection which
explains \eqn{symcaps} combinatorially.

\begin{example}
Let $m=5$ and $k=4$. 
There are seven $(5,4)$-capsid partitions of $16$:
\beqs
(1^{16}),\quad (4^{4}), \quad (1^{11},5), \quad (1^{6},5^2),
\quad (1,5^{3}), \quad (1^{6},10), \quad (1^2,14),
\eeqs
so that $a(5,4,16)=7$, thus (see Example~\ref{ex:51})
confirming $a(5,4,16)=a(5,1,16)$.                                                         
\end{example}

Apart from the concrete applications of $(m,k)$-capsids in this paper,
we believe that this new class of partitions is an object of interest
worthy of independent study from various (including combinatorial
and analytic) points of view.
In particular, it would be interesting to establish
asymptotic formulas for the number of $(m,k)$-capsid partitions
of $n$, as $n$ goes to infinity. For $t$-core partitions, this has
been established by Anderson~\cite{An}.

\section{Special $q$-Series identities involving congruences}
\label{sectspecial}

The two Rogers--Ramanujan functions
\beqs
G(q):=\sum_{k=0}^\infty\frac{q^{k^2}}{(q;q)_k}\qquad\text{and}\qquad
H(q):=\sum_{k=0}^\infty\frac{q^{k(k+1)}}{(q;q)_k}
\eeqs
satisfy a number of remarkable properties.
First of all, they factorize in closed form,
\beqs
G(q)=\frac1{(q,q^4;q^5)_\infty}\qquad\text{and}\qquad
H(q)=\frac1{(q^2,q^3;q^5)_\infty},
\eeqs
which are the two Rogers--Ramanujan identities~\cite{Ro}.

We will make essential use of the beautiful identity
\begin{equation}\label{gh55}
1=H(q)G(q^{11})-q^2G(q)H(q^{11}),
\end{equation}
which Ramanujan stated without proof in one of his letters to
Hardy~\cite{R2}. The first published proof is due to Rogers~\cite{Ro2}.
Ramanujan also gave another remarkable identity related to \eqn{gh55}:
\beq\label{gh55s}
1=H(q)G(q)^{11}-q^2G(q)H(q)^{11}-11qG(q)^6H(q)^6.
\eeq

Ramanujan's pair of equations \eqn{gh55} and \eqn{gh55s} can be
compared with the similar pair of identities
\beq\label{gh39}
1=\twidit{H}(q)\twidit{G}(q^3)-q^2\twidit{G}(q)\twidit{H}(q^3),
\eeq
and
\beqs
1=\twidit{H}(q)\twidit{G}(q)^3-q^2\twidit{G}(q)\twidit{H}(q)^3-
3q\twidit{G}(q)^2\twidit{H}(q)^2,
\eeqs
where
\begin{align*}
\twidit{G}(q)&=\frac1{(q,q^3,q^4,q^9,q^{10},q^{12};q^{13})_\infty},\\
\intertext{and}
\twidit{H}(q)&=\frac1{(q^2,q^5,q^6,q^7,q^8,q^{11};q^{13})_\infty},
\end{align*}
proved in S.~Robins' 1991 thesis~\cite{R}.

We define
\beq\label{eq:pmk}
{P}_{m,k}(q):= (q^k; q^m)_\infty (q^{m-k};q^m)_\infty.
\eeq
On one hand the symmetry
\beq\label{symp}
P_{m,k}(q)=P_{m,m-k}(q)
\eeq
and the scaling of parameters
\beqs
P_{sm,sk}(q)=P_{m,k}(q^s)
\eeqs
are obvious.
On the other hand, using
$(a;q)_\infty=(a,aq;q^2)_\infty$ and \eqn{symp}, we also have
\beq\label{doubling}
P_{m,k}(q)=P_{2m,k}(q)P_{2m,m-k}(q).
\eeq

Ramanujan's identity \eqn{gh55} can be written as
\beq\label{gh55p}
1=\frac{1}{{P}_{5,2}(q) {P}_{55,11}(q)} -
\frac{q^2}{{P}_{5,1}(q){P}_{55,22}(q)},
\eeq
or, in view of \eqn{doubling}, as
\beq\label{gh55pd}
1=\frac{1}{{P}_{10,2}(q){P}_{10,3}(q)
{P}_{110,11}(q){P}_{110,44}(q)} 
-
\frac{q^2}{{P}_{10,1}(q){P}_{10,4}(q)
{P}_{110,22}(q){P}_{110,33}(q)},
\eeq
which we will make use of.

Bressoud \cite{Br-thesis} gave an immediate combinatorial interpretation 
of \eqn{gh55p},
viewed as a mod $55$ identity. From a combinatorial perspective,
\eqn{gh55p} is a ``shifted partition identity''.
Other such identities (such as \eqn{gh39}), possibly related to other
moduli, have been discussed and established in \cite{A1,K,G2,GY}.
Most of the known shifted partition identities share the feature that they can
be written as some mod $m$ identity in the form
\begin{subequations}\label{AB}
\begin{equation}\label{ABid}
1=A_m(q)-q^dB_m(q)
\end{equation}
with
\beq
A_m(q)=\frac 1{\prod_{j=1}^{12}P_{m,k_j}(q)}\qquad\text{and}\qquad
B_m(q)=\frac 1{\prod_{j=1}^{12}P_{m,l_j}(q)},
\eeq
\end{subequations}
for some $m,d\in\N$ and $0< k_j,l_j<m$, for $j=1,\dots,12$.

It is clear that \eqn{gh55p} has the form of \eqn{AB},
obtained by taking $m=55$, $d=2$,
$k_j=5j-3$ and $l_j=5j-4$ for $j=1\dots,11$, $k_{12}=11$,
and $l_{12}=22$.
Similarly, \eqn{gh39} is a special case of \eqn{AB},
obtained by taking $m=39$, $d=2$,
$k_j=13j-11$ and $l_j=13j-12$ for $j=1,2,3$,
$k_j=13j-47$ and $l_j=13j-49$ for $j=4,5,6$,
$k_j=13j-85$ and $l_j=13j-87$ for $j=7,8,9$,
$k_{10}=3$, $k_{11}=9$, $k_{12}=12$,
$l_{10}=6$, $l_{11}=15$, and $l_{12}=18$.

\section{Combinatorial interpretations of the tau function}

We are ready to collect the ingredients we prepared
and apply them to obtain combinatorial
interpretions of Ramanujan's tau function.

First observe that by Proposition~\ref{prop:caps} and the notation
\eqn{eq:pmk}, the generating function for
$(m,k)$-capsid partitions takes the form
\beqs
\cpgn{m}{k}{q}=\frac{(q^m;q^m)_\infty}{P_{m,k}(q)}.
\eeqs
Now, we multiply both sides of \eqn{gh55pd} by 
$q^{110}(q^{110};q^{110})_\infty^{24}$
and rewrite the resulting identity in terms of $(10,k)$-capsid and
$11$-core generating functions:
\begin{align}\label{eq:tauid}
\sum_{n=1}^\infty \tau(n) q^{110n}
&=q^{110}\,\frac{(q^{10};q^{10})_\infty^2}{{P}_{10,2}(q){P}_{10,3}(q)}
\frac{(q^{110};q^{110})_\infty^2}{{P}_{10,1}(q^{11}){P}_{10,4}(q^{11})}
\frac{(q^{110};q^{110})_\infty^{22}}{(q^{10};q^{10})_\infty^2}\notag\\
&\quad-
q^{112}\,\frac{(q^{10};q^{10})_\infty^2}{{P}_{10,1}(q){P}_{10,4}(q)}
\frac{(q^{110};q^{110})_\infty^2}{{P}_{10,2}(q^{11}){P}_{10,3}(q^{11})}
\frac{(q^{110};q^{110})_\infty^{22}}{(q^{10};q^{10})_\infty^2}
\notag\\
&=q^{110}\,\cpgn{10}{2}{q}\,\cpgn{10}{3}{q}\,\cpgn{10}{1}{q^{11}}\,
\cpgn{10}{4}{q^{11}}\,
\tcgn{11}{q^{10}}^2
\notag \\
&\quad-
q^{112}\,\cpgn{10}{1}{q}\,\cpgn{10}{4}{q}\,\cpgn{10}{2}{q^{11}}\,
\cpgn{10}{3}{q^{11}}\,
\tcgn{11}{q^{10}}^2 .
\end{align}

We let $\mathcal{C}_{m,k}$ be the set of $(m,k)$-capsids
and $\mathcal{T}_t$ be the set of $t$-cores.
We define two sets of vector partitions:
\begin{align*}
\mathcal{A} &:= 
\cps{10}{2}\times \cps{10}{3}\times \cps{10}{1}\times \cps{10}{4}\times 
\tcs{11} \times \tcs{11},\\
\mathcal{B} &:= 
\cps{10}{1}\times \cps{10}{4}\times \cps{10}{2}\times \cps{10}{3}\times 
\tcs{11} \times \tcs{11}.
\end{align*}
For a partition $\pi$ we let $\abs{\pi}$ denote the sum of parts.
For $\vec{\pi} = (\pi_1,\pi_2,\pi_3,\pi_4,\pi_5,\pi_6)$ in $\mathcal{A}$
or in  $\mathcal{B}$ we define
\beqs
\abs{\vec{\pi}}:= \abs{\pi_1} + \abs{\pi_2} + 11 \cdot \abs{\pi_3} 
+ 11\cdot \abs{\pi_4} + 10\cdot\abs{\pi_5} + 10\cdot\abs{\pi_6}.
\eeqs
If $\abs{\vec{\pi}}=n$ we say that $\vec{\pi}$ is a vector partition of $n$.
We let
\beqs
a(n) := \mbox{the number of vector partitions in $\mathcal{A}$ of $n$},
\eeqs
and
\beqs
b(n) := \mbox{the number of vector partitions in $\mathcal{B}$ of $n$}.
\eeqs
Then it is clear
that
\beqs
1 + \sum_{n=1}^\infty a(n) q^n =
\cpgn{10}{2}{q}\,\cpgn{10}{3}{q}\,\cpgn{10}{1}{q^{11}}\,\cpgn{10}{4}{q^{11}}\,
\tcgn{11}{q^{10}}^2,
\eeqs
and
\beqs
1 + \sum_{n=1}^\infty b(n) q^n =
\cpgn{10}{1}{q}\,\cpgn{10}{4}{q}\,\cpgn{10}{2}{q^{11}}\,\cpgn{10}{3}{q^{11}}\,
\tcgn{11}{q^{10}}^2.
\eeqs

We define $a(0)=1$, $b(-2)=b(-1)=0$, and $b(0)=1$.
From \eqn{eq:tauid} we readily obtain the following result:
\begin{theorem}[A combinatorial formula for Ramanujan's tau function]
\label{thm1}
For $n\ge1$, we have
\begin{alignat*}{2}
(i)&\qquad&
\tau(n)&= a(110n - 110) - b(110n - 112),\\
\intertext{and}
(ii)&&
0&=a(n) - b(n-2),\qquad\qquad
\text{if $n\not\equiv0$ mod $110$}.
\end{alignat*}
\end{theorem}

This formula for $\tau(n)$ has the advantage that it is compact and
combinatorial,
but the numbers obtained when
computing $a(m)$ and $b(m)$ get large.
\begin{example}\label{extau2}
\beqs
\tau(2)=a(110)-b(108)=174780-174804=-24.
\eeqs
\end{example}

Let us look for another formula for $\tau(n)$ using a smaller modulus
than $110$, say $10$. For this we use Ramanujan's
\eqn{gh55s} which we rewrite in the form
\begin{align*}
1=\frac 1{P_{10,2}(q)P_{10,3}(q)P_{10,1}(q)^{11}P_{10,4}(q)^{11}}
-{}\frac {q^2}{P_{10,1}(q)P_{10,4}(q)P_{10,2}(q)^{11}P_{10,3}(q)^{11}}&\notag\\
{}-{}\frac {11\,q}{P_{10,1}(q)^6P_{10,4}(q)^6P_{10,2}(q)^6P_{10,3}(q)^6}&.
\end{align*}
We multiply both sides of this relation by 
$q^{10}(q^{10};q^{10})_\infty^{24}$
and rewrite the resulting identity in terms of $(10,k)$-capsid
generating functions:
\begin{align}\label{eq:tauid2}
\sum_{n=1}^\infty \tau(n) q^{10n}
&=q^{10}\,\cpgn{10}{2}{q}\,\cpgn{10}{3}{q}\,\cpgn{10}{1}{q}^{11}\,
\cpgn{10}{4}{q}^{11}
\notag\\
&\quad-
q^{12}\,\cpgn{10}{1}{q}\,\cpgn{10}{4}{q}\,\cpgn{10}{2}{q}^{11}\,
\cpgn{10}{3}{q}^{11}\notag\\[.5em]
&\quad-
11\,q^{11}\,\cpgn{10}{1}{q}^6\,\cpgn{10}{4}{q}^6\,\cpgn{10}{2}{q}^6\,
\cpgn{10}{3}{q}^6.
\end{align}

We now define three sets of vector partitions, each defined by a
twentyfour-fold Cartesian product:
\beqs
\mathcal{U} :=\Cross_{j=1}^{24}\;\mathcal{U}_j,\qquad
\mathcal{V} :=\Cross_{j=1}^{24}\;\mathcal{V}_j,\qquad
\mathcal{W} :=\Cross_{j=1}^{24}\;\mathcal{W}_j,
\eeqs
where
\begin{alignat*}4
\mathcal{U}_1&=\cps{10}{2},&\quad \mathcal{U}_2&=\cps{10}{3},&\quad
\mathcal{U}_j&=\cps{10}{1},&\quad \mathcal{U}_{j+11}&=\cps{10}{4},\\
\mathcal{V}_1&=\cps{10}{1},&\quad \mathcal{V}_2&=\cps{10}{4},&\quad
\mathcal{V}_j&=\cps{10}{2},&\quad \mathcal{V}_{j+11}&=\cps{10}{3},
\quad\text{for $j=3,\dots,13$},\\
\mathcal{W}_l&=\cps{10}{1},&\quad \mathcal{W}_{l+6}&=\cps{10}{4},&\quad
\mathcal{W}_{l+12}&=\cps{10}{2},&\quad \mathcal{W}_{l+18}&=\cps{10}{3},
\quad\text{for $l=1,\dots,6$}.
\end{alignat*}
For $\vec{\pi} = (\pi_1,\pi_2,\dots,\pi_{24})$ in $\mathcal{U}$,
$\mathcal{V}$, or in $\mathcal{W}$ we define
\beqs
\abs{\vec{\pi}}:= \sum_{j=1}^{24}\abs{\pi_j}.
\eeqs
If $\abs{\vec{\pi}}=n$ we say that $\vec{\pi}$ is a vector partition of $n$.
We let
\begin{align*}
u(n)&:= \mbox{the number of vector partitions in $\mathcal{U}$ of $n$},\\
v(n)&:= \mbox{the number of vector partitions in $\mathcal{V}$ of $n$},\\
w(n)&:= \mbox{the number of vector partitions in $\mathcal{W}$ of $n$}.
\end{align*}
Then it is clear
that
\begin{align*}
1 + \sum_{n=1}^\infty u(n) q^n &=
\cpgn{10}{2}{q}\,\cpgn{10}{3}{q}\,\cpgn{10}{1}{q}^{11}\,\cpgn{10}{4}{q}^{11},\\
1 + \sum_{n=1}^\infty v(n) q^n &=
\cpgn{10}{1}{q}\,\cpgn{10}{4}{q}\,\cpgn{10}{2}{q}^{11}\,\cpgn{10}{3}{q}^{11},\\
1 + \sum_{n=1}^\infty w(n) q^n &=
\cpgn{10}{1}{q}^6\,\cpgn{10}{4}{q}^6\,
\cpgn{10}{2}{q}^{6}\,\cpgn{10}{3}{q}^{6}.
\end{align*}

We shall define $u(0)=1$, $v(-2)=v(-1)=0$, $v(0)=1$, $w(-1)=0$ and $w(0)=1$.
From \eqn{eq:tauid2} we readily obtain the following result:
\begin{theorem}[Another combinatorial formula for Ramanujan's tau function]
\label{thm2}
 For $n\ge1$, we have
\begin{alignat*}{2}
(i)&\qquad&
\tau(n)&= u(10n - 10) - v(10n - 12)-11w(10n-11),\\
\intertext{and}
(ii)&&
0&=u(n) - v(n-2)-11w(n-1),\qquad\qquad
\text{if $n\not\equiv0$ mod $10$}.
\end{alignat*}
\end{theorem}

Computationally, this formula is not better than the one
in Theorem~\ref{thm1}, as the numbers obtained when
computing $u(l)$, $v(l)$, and $w(l)$ get large as well.
In particular (compare with Example~\ref{extau2})
\beqs
\tau(2)=u(10)-v(8)-11w(9)=381405-3139-11\cdot 34390=-24.
\eeqs

\section{Other combinatorial interpretations}

It is clear that one can also employ other shifted partition
identities (which are listed in \cite{GY}) to obtain
combinatorial interpretations of Ramanujan's tau function
which are similar to those in Theorems~\ref{thm1} and \ref{thm2}.
Also, one can in principle use those shifted partition
identities to obtain combinatorial interpretations
for powers other than $24$ of Dedekind's eta function.
On the contrary, our method can also be used to establish
that certain sets of Cartesian products 
involving (restricted) partitions, $(m,k)$-capsid partitions and
$t$-core partitions have the same cardinality.

For example, Jacobi showed that
\beq
(q;q)_\infty^3=\sum_{n=0}^\infty(-1)^n(2n+1)q^{13\binom{n+1}2}.
\label{jacid}
\eeq
Joichi and Stanton \cite{Jo-St} have found a combinatorial proof of
\eqn{jacid}.
Combining the $q\mapsto q^{13}$ case of this identity with \eqn{gh39}, 
we have
\beqs
\sum_{n=0}^\infty(-1)^n(2n+1)q^{13\binom{n+1}2}=
\frac{\cpgn{13}{2}{q}\cpgn{13}{5}{q}\cpgn{13}{6}{q}}
{P_{13,1}(q^3)P_{13,3}(q^3)P_{13,4}(q^3)}-q^2
\frac{\cpgn{13}{1}{q}\cpgn{13}{3}{q}\cpgn{13}{4}{q}}
{P_{13,2}(q^3)P_{13,5}(q^3)P_{13,6}(q^3)},
\eeqs
which leads to the following combinatorial interpretation:

Let $\mathcal{D}$ and $\mathcal{E}$ be two sets of vector partitions,
defined by
\begin{align*}
\mathcal{D} &:= 
\cps{13}{2}\times \cps{13}{5}\times \cps{13}{6}\times X,\\
\mathcal{E} &:= 
\cps{13}{1}\times \cps{13}{3}\times \cps{13}{4}\times Y,
\end{align*}
where $X$ is the set of partitions with parts
$\pm 1$, $\pm 3$, $\pm 4$ mod $13$, and
$Y$ is the set of partitions with parts
$\pm 2$, $\pm 5$, $\pm 6$ mod $13$.
For $\vec{\pi} = (\pi_1,\pi_2,\pi_3,\pi_4)$ in $\mathcal{D}$
or in  $\mathcal{E}$ we define
\beqs
\abs{\vec{\pi}}:= \abs{\pi_1} + \abs{\pi_2} + \abs{\pi_3} +
3 \cdot \abs{\pi_4}.
\eeqs
If $\abs{\vec{\pi}}=n$ we say that $\vec{\pi}$ is a vector partition of $n$.
We let
\beqs
d(n) := \mbox{the number of vector partitions in $\mathcal{D}$ of $n$},
\eeqs
and
\beqs
e(n) := \mbox{the number of vector partitions in $\mathcal{E}$ of $n$}.
\eeqs
Then, with $e(-1)=0$ and $e(0)=1$,
we have the following result involving triangular numbers:
\begin{theorem}
\label{thm3}
For $m\ge1$, we have
\begin{alignat*}{3}
(i)&\qquad&
d(m)-e(m-2)&=(-1)^n(2n+1),
\qquad&&\text{if $m=\frac{13n(n+1)}2$}\\
\intertext{and}
(ii)&&
d(m)-e(m-2)&=0,\qquad
&&\text{otherwise}.
\end{alignat*}
\end{theorem}

\section{Capsid symmetries} 
\label{sectsyms}
In this section we explain 
\eqn{symcaps} combinatorially.
In Section \ref{sectptns} we defined $(m,k)$-capsids.
Now we define more general capsids. Suppose $m$, $r_1$, $r_2\in\N$,
$m\ge 2$ and $0 < r_1 \ne r_2 < m$. We say that a partition $\pi$
is an $(m,r_1,r_2)$-capsid if, and only if, the possible
parts of $\pi$ are $r_1$ or are congruent to $0$ or $r_2$ 
mod $m$, and satisfy the following two conditions:
\begin{enumerate}
\item[(i)]
if $\mu_{r_1}=0$, i.e., $r_1$ is not a part,
then all parts are congruent to $r_2$ mod $m$;
\item[(ii)]
If $\mu_{r_1}>0$, then $r_1$ is the smallest part
and the largest part congruent to $0$ mod $m$ is
$\le m\cdot\mu_{r_1}$ and all parts congruent to $r_2$ mod $m$
are $> m\cdot\mu_{r_1}$.
\end{enumerate}

Thus (for $k\ne m/2$) a $(m,k)$-capsid is a $(m,m-k,k)$-capsid.

\medskip
Let $\gamma(m,r_1,r_2,n)$ be the number of $(m,r_1,r_2)$-capsid 
partitions of $n$.
We have the
following generating function
\begin{equation}
{C}_{m,r_1,r_2} (q)= 1 + \sum_{n=1}^\infty \gamma(m,r_1,r_2,n) q^n
=
\sum_{n=0}^\infty \frac{q^{r_1 n}}
{(q^m;q^m)_n (q^{mn+r_2};q^m)_\infty}.
\end{equation}

For an $(m,r_1,r_2)$-capsid $\pi$ we define two statistics:
\begin{enumerate}
\item[(i)]
We let $\alpha(\pi)=\mu_{r_1}$.
\item[(ii)]
We let $\beta(\pi)$ be the number of parts of $\pi$ that
are congruent to $r_2$ mod $m$.
\end{enumerate}

Let $\gamma(m,r_1,r_2,a,b,n)$ be the number of 
$(m,r_1,r_2)$-capsid partitions $\lambda$ of $n$ with 
$\alpha(\lambda)=a$ and $\beta(\lambda)=b$.

Let $\mathcal{C}_{m,r_1,r_2}$ be the set of $(m,r_1,r_2)$-capsids.
It is clear that we have the following generating function

\begin{equation}\label{gencapsid}
{C}_{m,r_1,r_2} (x,y,q)= 
\sum_{\lambda\in\mathcal{C}_{m,r_1,r_2}} x^{\alpha(\lambda)}
y^{\beta(\lambda)} q^{\abs{\lambda}}
=
\sum_{n=0}^\infty \frac{x^n q^{r_1n}}
{(q^m;q^m)_n (yq^{mn+r_2};q^m)_\infty}.
\end{equation}
By \eqn{Cgenf} we have the following result:
\begin{proposition}\label{prop:gencaps}
The generating function for $(m,r_1,r_2)$-capsids is
\beqs
{C}_{m,r_1,r_2}(x,y,q)=\frac{ (xyq^{r_1+r_2};q^m)_\infty }
{ (xq^{r_1};q^m)_\infty (yq^{r_2};q^m)_\infty}.
\eeqs
\end{proposition}

Since $C_{m,r_1,r_2}(x,y,q) = C_{m,r_2,r_1}(y,x,q)$ and
$$
\gamma(m,r_1,r_2,a,b,n) = \mbox{the coefficient of $x^a y^b q^n$
 in $C_{m,r_1,r_2}(x,y,q)$},
$$
we deduce that
\beq
\label{gensymcaps}
\gamma(m,r_1,r_2,a,b,n) =
\gamma(m,r_2,r_1,b,a,n) 
\eeq
for all $n$. By summing over $a$, $b$ this implies that
\beq
\label{gensymcaps2}
\gamma(m,r_1,r_2,n) = \gamma(m,r_2,r_1,n) 
\eeq
for all $n$, which generalizes \eqn{symcaps}.

In this section we give a combinatorial proof of \eqn{gensymcaps}
by constructing a bijection
$$
\mathcal{J}\,:\,\mathcal{C}(m,r_1,r_2,a,b,n)
\longrightarrow
\mathcal{C}(m,r_2,r_1,b,a,n),
$$
where
$$
\mathcal{C}(m,r_1,r_2,a,b,n):=\{
\lambda\in\mathcal{C}_{m,r_1,r_2}\,:\,
\alpha(\lambda)=a,\,
\beta(\lambda)=b,\,    
\abs{\lambda}=n\}.
$$
For a given $\lambda\in\mathcal{C}(m,r_1,r_2,a,b,n)$ 
we divide each part
of $\lambda$ congruent to $0$ mod $m$ by $m$ to form a partition
$\pi_1$, and form a partition $\pi_2$ by subtracting $r_2$ from
 each part of $\lambda$ congruent to $r_2$ mod $m$ and dividing by
$m$. Then we form the diagram given below in Figure 1.

\begin{figure}[ht]
\centering
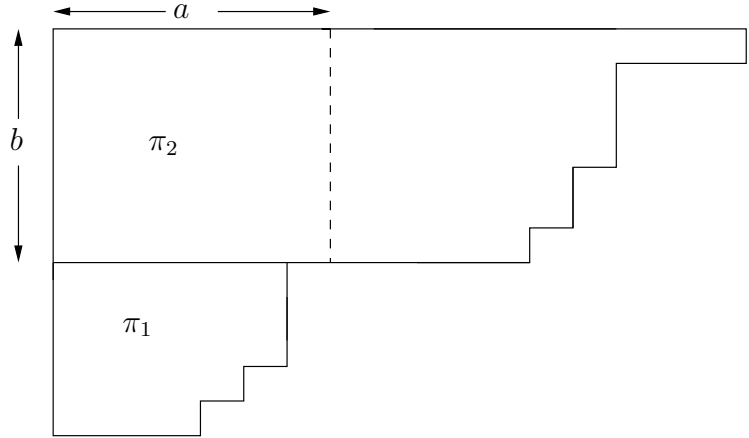
\caption{The partition $\pi = \pi_2 \cup \pi_1$}
\label{FIG1}
\end{figure}

\noindent
\textbf{Note:}
\begin{enumerate}
\item[(i)]
The partition $\pi_2$ has $b$ parts each $\ge a$.
\item[(ii)]
The parts of $\pi_1$ are $\le a$.
\item[(iii)]
$a$ is the number of smallest parts of $\lambda$.
\item[(iv)]
If $b=0$ then $\pi_2=(\,)$ is empty.
\item[(v)]
It is possible that $\pi_1=(\,)$ is empty.
\end{enumerate}
Next we form the partition
$$
\pi = \pi_2 \cup \pi_1,
$$
and take the conjugate $\pi'$. Then we form partitions
$\twidit{\pi}_1$, $\twidit{\pi}_2$, where $\pi' = \twidit{\pi_2}
\cup \twidit{\pi}_1$ and given in the diagram below in Figure 2.

\begin{figure}[ht]
\centering
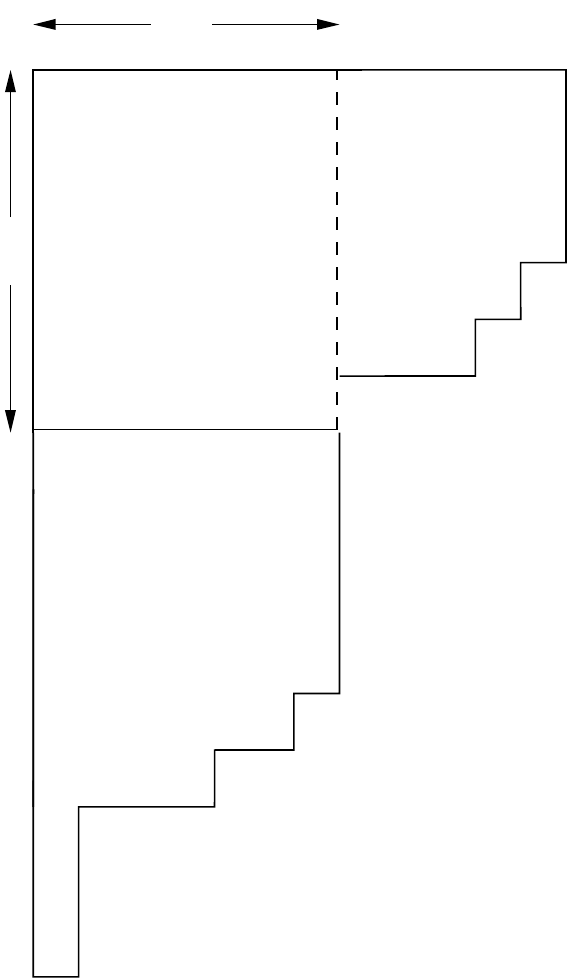
\caption{The partition $\pi' = \twidit{\pi}_2 \cup \twidit{\pi}_1$}
\label{FIG2}
\end{figure}

We form the partition $\twidit{\lambda}$ with $b$ smallest
parts $r_2$, with parts from $\twidit{\pi}_1$ multiplied
by $m$, and with adding $r_1$ to each part of $\twidit{\pi}_2$
after multiplying by $m$. We see that
$\twidit{\lambda}\in \mathcal{C}(m,r_2,r_1,b,a,n)$.
We define $\mathcal{J}(\lambda)=\twidit{\lambda}$ and
$\mathcal{J}$ is a bijection.

\vfill
\eject

\begin{example}
\label{ex:bijeg}
$\lambda=(1^3,5^1,15^2,22^1,27)$, $m=5$, $r_1=1$, $r_2=2$, $a=3$,
$b=2$ and $n=87$. From Figure  3 below we see that
$\mathcal{J}(\lambda)=\twidit{\lambda}
=(2^2,5^1, 10^1, 21^2, 26)\in\mathcal{C}(5,2,1,2,3,87)$.
\end{example}

\begin{figure}[ht]
\centering
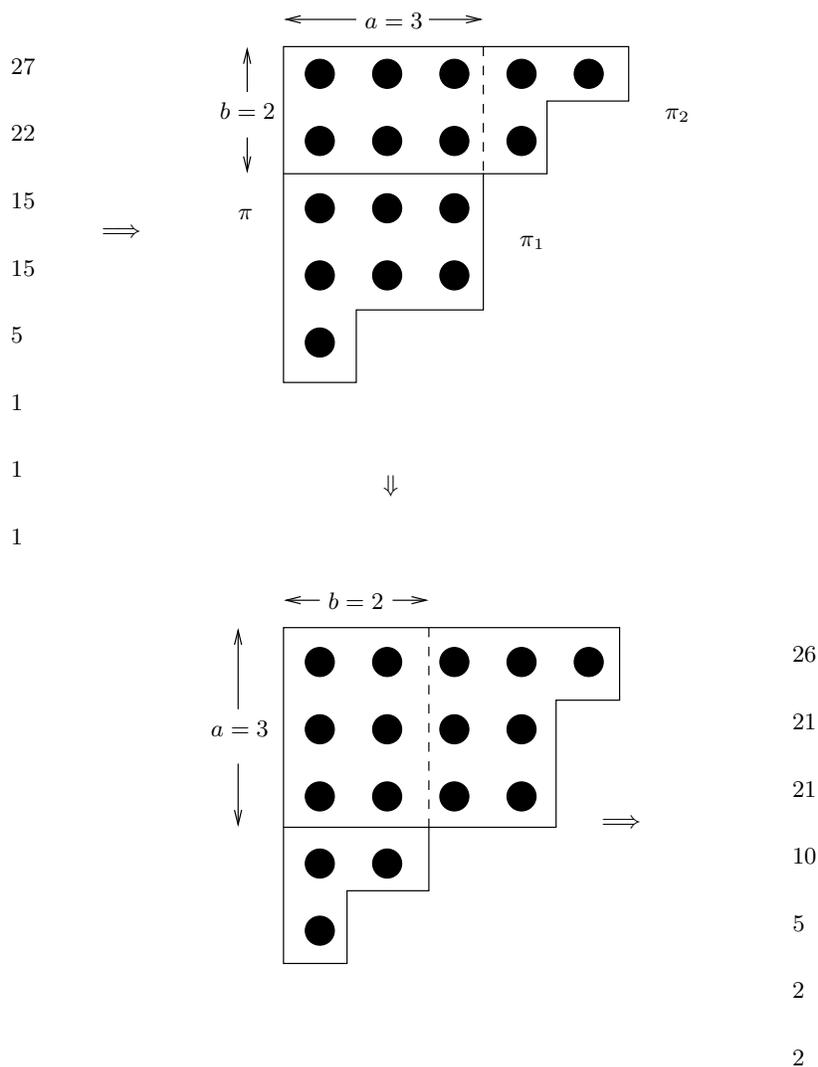
\caption{Example of bijection $\mathcal{J}$}
\label{FIG3}
\end{figure}


\section{Concluding remarks}
\label{sec:conclusion}

The main goal of this paper was to give a combinatorial interpretation
of the Ramanujan tau function $\tau(n)$ in terms of partitions.
In Theorems \ref{thm1} and \ref{thm2} we were able to write $\tau(n)$ as the
difference of two, resp.\ three  partition functions with 
nonnegative coefficients.
We do not claim that this should be used to calculate $\tau(n)$ or to find
improved asymptotics. 
Charles \cite{Ch} and Edixhoven \cite{Ed} have found
efficient algorithms for computing $\tau(n)$.
We note that there are other ways to express $\tau(n)$ in terms
of nonnegative arithmetic functions.
For example it is well known that there are many identities for $\tau(n)$ 
in terms of divisor functions.
Following Ramanujan \cite[p.163]{R1} we define
$$
\Phi_j(q) = \sum_{n=1}^\infty \frac{n^j q^n}{1-q^n} = 
\sum_{m,n\ge1} n^j q^{nm} = \sum_{n=1}^\infty \sigma_j(n) q^n,
$$
for $j\ge 1$ odd and where $\sigma_j(n) = \sum_{d\mid n} d^j$.
For $n\ge2$ even we define the Eisenstein series
\beq
E_n(q) 
= 1 - \frac{2n}{B_n} \Phi_{n-1}(q),
\label{eq:En}
\eeq
where $B_n$ is the $n$-th
Bernoulli number.
Then it is well-known that
$$
\sum_{n=1}^\infty \tau(n) q^n = \frac{1}{1728}\left( E_4(q)^3 - E_6(q)^2
\right).
$$
See \cite[p.172]{Mi01} for a history
of this identity.
Thus $\tau(n)$ can be written in terms of the divisor sums $\sigma_3(n)$
and $\sigma_5(n)$.  A nicer identity is  
$$
\sum_{n=1}^\infty \tau(n) q^n = \frac{691}{762048}\left( E_{12}(q) - E_6(q)^2
\right),
$$
which implies that
\beq
\tau(n) = \frac{65}{756} \sigma_{11}(n) 
 + \frac{691}{756} \sigma_5(n) 
 - \frac{691}{3} \sum_{m=1}^{n-1} \sigma_5(m) \sigma_5(n-m).
\label{tau691}
\eeq
See \cite[Eq.(1.15)]{Mi02}.
It is not clear how to interpret the right side of \eqn{tau691}
combinatorially, although
it does express $\tau(n)$ in terms of nonnegative arithmetic functions.


\end{document}